\documentclass{article}
\usepackage{amsfonts}
\usepackage{amsmath}
\usepackage{amssymb}

\newtheorem{theorem}{Theorem}

\newtheorem{lemma}[theorem]{Lemma} 
\newtheorem{corollary}[theorem]{Corollary}

\newtheorem{remark}[theorem]{Remark}

\newcommand{\barS}{{\overline{S}}}
\newcommand{\barm}{{\overline{m}}}

\newcommand{\Lambdahat}{{\widehat\Lambda}}
\newcommand{\lhat}{{\hat\ell}}

\newcommand{\zz}{{\mathbb Z}}
\newcommand{\gen}[1]{\left\langle #1\right\rangle}
\newcommand{\set}[1]{\left\{#1\right\}}

\title{Artin HNN-extensions virtually embed in Artin groups}

\author{T. Hsu\thanks{This work was done at MSRI,
    Berkeley, where research is supported in part by NSF grant DMS-0441170}
    \and I. J. Leary$^*$\thanks{partially
      supported by NSF grant DMS-0505471.}}

\date{\today}

\newenvironment{proof}[1][]{\begin{trivlist} \item[\hskip\labelsep
\emph{Proof#1.}]}{\foorp \end{trivlist}}
\newcommand{\foorp}{{\unskip\nobreak\hfil\penalty50
 \hskip1em\vadjust{}\nobreak\hfil \vrule height3pt width3pt depth0pt
 \parfillskip=0pt \finalhyphendemerits=0 \par}}

\begin{document} 

\maketitle

\begin{abstract} 
  An Artin HNN-extension is an HNN-extension of an Artin group in
  which the stable letter conjugates a pair of suitably chosen subsets
  of the standard generating set.  We show that some finite index
  subgroup of an Artin HNN-extension embeds in an Artin group.  We
  also obtain an analogous result for Coxeter groups.
\end{abstract}



\section{Results}
\label{results}

Recall that an {\sl Artin system} is a finite set $S$, together with a
function $m$ from the 2-element subsets of $S$ to the set
$\{2,3,4,\ldots\}\cup\{\infty\}$.  It may be helpful to consider $m$
as an edge-labelling on the complete graph with vertex set $S$.  The
{\sl Artin group} $A=A(S,m)$ corresponding to an Artin system $(S,m)$
is defined by a presentation.  The generators for $A$ are the elements
of the set $S$.  There is one relation in the presentation for each
unordered pair $s,t\in S$ such that $m(t,s)=m(s,t)<\infty$, of the
form
\begin{equation}\label{eq:artin-relation}
  sts\cdots = tst\cdots,
\end{equation}
where there are $m(s,t)$ symbols on each side of
\eqref{eq:artin-relation}.  Call these relations the {\sl Artin
  relations} associated to $(S,m)$.  If $S'$ is a subset of $S$, and
we write $m$ for the restriction to $S'$ of $m$, inclusion induces a
group homomorphism from $A(S',m)$ to $A(S,m)$, and this homomorphism
is known to be injective \cite{vdl,paris}.  The image of such a
homomorphism is called an {\sl Artin subgroup} of $A(S,m)$.

An Artin system also gives rise to a {\sl Coxeter group} $W=W(S,m)$.
The generators of the Coxeter group are the elements of $S$, subject
to the Artin relations for $(S,m)$ together with the relations $s^2=1$
for all $s\in S$.  Call these relations the {\sl Coxeter relations}
for $W$.  Note that in the presence of the relations $s^2=1=t^2$, the
Artin relation between $s$ and $t$ is equivalent to the relation
$(st)^{m(s,t)}=1$.  For $S'$ a subset of $S$, the natural map
$W(S',m)\rightarrow W(S,m)$ is injective, and the image of such a
homomorphism is called a {\sl Coxeter subgroup} of $W(S,m)$
\cite{bourbaki,brown}.

The {\sl set of labels}, $L(S,m)$, for an Artin system $(S,m)$ is the
image of $m$, viewed as a subset of
$\set{2,3,4,\ldots}\cup\set{\infty}$.  An Artin system $(S,m)$ is
called right-angled if $L(S,m)\subseteq \{2,\infty\}$, and $A(S,m)$
(resp.\ $W(S,m)$) is then called a {\sl right-angled Artin group}
(resp.\ {\sl right-angled Coxeter group}).

Now suppose $\phi:S'\rightarrow S''$ is a bijection between two
subsets of $S$ that preserves labels in the sense that
\begin{equation}\label{eq:label-preserving}
  m(s,s')=m(\phi(s),\phi(s')) \qquad\qquad\text{for all $s,s'\in S'$.}
\end{equation}
The {\sl Artin HNN-extension} $G(S,m,\phi)$ (resp.\ {\sl Coxeter
  HNN-extension} $C(S,m,\phi)$) associated to $(S,m,\phi)$ is the free
product of $A(S,m)$ (resp.\ $W(S,m)$) and an infinite cyclic group
$\gen{t}$, modulo the relations $t^{-1}st=\phi(s)$ for all $s\in S'$.
Thus an Artin HNN-extension (resp.\ a Coxeter HNN-extension) is a
special type of HNN-extension in which the base group is an Artin
(resp.\ a Coxeter) group and the stable letter conjugates an Artin
(resp.\ a Coxeter) subgroup to another Artin (resp.\ Coxeter) subgroup
via the map induced by a bijection of standard generating sets.  The
{\sl set of labels} for $G(S,m,\phi)$ (resp.\ $C(S,m,\phi)$) is the
set of labels for the group $A(S,m)$ (resp.\ $W(S,m)$).

We require one more definition before stating our main result.  For
any property $\cal P$ of groups, a group $G$ {\sl virtually has
  property $\cal P$} if there exists a finite-index subgroup $H\leq G$
such that $H$ has property $\cal P$.

\begin{theorem}
  Any Artin (resp.\ Coxeter) HNN-extension virtually embeds in an
  Artin (resp.\ Coxeter) group.  If $L$ is the set of labels for the
  HNN-extension, then the Artin (resp.\ Coxeter) group in which it
  virtually embeds can be chosen to have its set of labels equal to
  $L\cup \{2,\infty\}$.  In particular right-angled HNN-extensions
  virtually embed in right-angled groups.
  \label{main} 
\end{theorem}

Theorem~\ref{main} yields the following application, which was our
original motivation for considering Artin HNN-extensions.

\begin{corollary}\label{counterex}
  For any finite group $Q$ not of prime power order, there exists $n$
  and a subgroup $G$ of $SL(n,\zz)$ such that $G$ virtually has a
  finite classifying space and $G$ has infinitely many conjugacy 
  classes of subgroups isomorphic to $Q$.
\end{corollary} 

\begin{proof} Note that by `$G$ virtually has a finite classifying 
  space' we mean that there is a finite index subgroup $H\leq G$ which
  has a finite $K(H,1)$.  In \cite{leary}, a group $G$ is constructed
  which virtually has a finite classifying space and which contains
  infinitely many conjugacy classes of subgroup isomorphic to $Q$.  By
  \cite[Thm.\ 28]{leary}, we know that $G$ virtually embeds in a
  right-angled Artin HNN-extension.  But then, by Theorem~\ref{main},
  $G$ virtually embeds in a right-angled Artin group, and
  by~\cite{dj,hsuwise}, any right-angled Artin group embeds in
  $SL(m,\zz)$ for some $m$, so the group $G$ embeds in $SL(n,\zz)$ for
  some (possibly larger)~$n$.
\end{proof} 

\begin{remark}
  Corollary \ref{counterex} was proved in \cite{learynucinkis} in the
  special case when $Q$ is not of the form
  ($p$-group)-by-cyclic-by-($q$-group) for primes $p$ and $q$.
\end{remark}

\section{Proofs}

To review the standard definitions from \cite[I.5]{trees}, recall that
a {\sl graph of groups} $(G_{(-)},\Gamma)$ consists of a directed
CW-graph $\Gamma$, together with a group $G_v$ for each vertex $v$ of
$\Gamma$ (the {\sl vertex groups}), a group $G_e$ for each edge $e$ of
$\Gamma$ (the {\sl edge groups}), and for each edge $e$, two injective
group homomorphisms $\iota_e:G_e\rightarrow G_{\iota(e)}$ and
$\tau_e:G_e\rightarrow G_{\tau(e)}$ (the {\sl attaching
  homomorphisms}), where $\iota(e)$ denotes the initial vertex of $e$
and $\tau(e)$ denotes the terminal vertex of $e$.  We define
$F(G_{(-)},\Gamma)$ to be the free product of the $G_v$ and infinite
cyclic groups $\gen{t_e}$, where $e$ runs over all edges of $\Gamma$,
modulo the relations
\begin{align}\label{eq:F}
  t_e^{-1}\iota_e(g) t_e &= \tau_e(g) &\text{for all $e$, $g\in G_e$;}
\end{align}
and for a basepoint $P_0$ in $\Gamma$, we define
$\pi_1(G_{(-)},\Gamma,P_0)$ to be the subgroup of all words in
$F(G_{(-)},\Gamma)$ that, after ignoring all syllables from $G_v$'s,
define paths in $\Gamma$ that begin and end at $P_0$.  Similarly, for
a spanning tree $T$ of $\Gamma$, we define $\pi_1(G_{(-)},\Gamma,T)$
to be the free product of the $G_v$ and infinite cyclic groups
$\gen{t_e}$, where $e$ runs over all edges of $\Gamma$ not in $T$,
modulo the relations
\begin{align}
  \iota_e(g) &= \tau_e(g) &\text{for all $e\in T$,
    $g\in G_e$, and} \label{eq:inT} \\
  t_e^{-1}\iota_e(g) t_e &= \tau_e(g) &\text{for all $e\notin T$,
    $g\in G_e$.} \label{eq:notinT}
\end{align}
Up to isomorphism, $\pi_1(G_{(-)},\Gamma,P_0)$ is independent of $P_0$
and $\pi_1(G_{(-)},\Gamma,T)$ is independent of $T$, and the two
groups are isomorphic, so we call their isomorphism class the {\sl
  fundamental group of $(G_{(-)},\Gamma)$}, denoted by
$\pi_1(G_{(-)},\Gamma)$.  We also call the elements $t_e$ {\sl stable
  letters}.

Let $k$ be a positive integer, let $V$, $E$, and $V_*$ be groups, and
let $f:E\rightarrow V$, $g:E\rightarrow V$, and $h_i:V\rightarrow V_*$
($i\in\zz/k$) be injective group homomorphisms.
\begin{enumerate}
\item\label{HNN-graph} Let $\Delta_1$ be a self-loop.  Taking vertex
  group $V$, edge group $E$, and attaching maps $f$ and $g$, we get a
  graph of groups $(\set{V,E},\Delta_1)$ whose fundamental group is an
  HNN extension with base group $V$ and associated subgroup $E$.  Let
  $t$ be the associated stable letter.
\item\label{index-k} Let $\Delta_k$ be a (directed) $k$-cycle, with
  vertex set $\zz/k$ and edges $(i,i+1)$ ($i\in\zz/k$).  Taking the
  vertex (resp.\ edge) groups $V_i$ (resp.\ $E_i$) to be copies of $V$
  (resp.\ $E$), and the attaching homomorphisms $f_i:E_i\rightarrow
  V_i$ and $g_i:E_i\rightarrow V_{i+1}$ to be copies of $f$ and $g$,
  with all indices running over all $i\in\zz/k$, we get a graph of
  groups $(\set{V_i,E_i},\Delta_k)$.  Let $\set{t_i}$ be the set of
  stable letters in the group $F(\set{V_i,E_i},\Delta_k)$.
\item\label{ambient-group} Let $\Gamma_k$ be a $k$-leaved rose, with
  edges labelled by $\zz/k$, and let $V_i$, $E_i$, $f_i$, and $g_i$ be
  defined as in \ref{index-k}.  Taking vertex group $V_*$, edge groups
  the $E_i$, and attaching homomorphisms $h_i\circ f_i$ and
  $h_{i+1}\circ g_i$, we get a graph of groups
  $(\set{V_*,E_i},\Gamma_k)$.
\end{enumerate}

\begin{lemma}\label{graphlem}
  Let $K$ be the kernel of the homomorphism $\rho:
  \pi_1(\set{V,E},\Delta_1)\rightarrow\zz/k$ defined by $\rho(V)=0$
  and $\rho(t)=1$.  Then $K=\pi_1(\set{V_i,E_i},\Delta_k)$, and $K$
  embeds in $\pi_1(\set{V_*,E_i},\Gamma_k)$.
\end{lemma}

\begin{proof}
  Let $\theta:
  \pi_1(\set{V_i,E_i},\Delta_k,0)\rightarrow\pi_1(\set{V,E},\Delta_1)$
  be the homomorphism defined by $\theta(v)=v\in V$ for all $v\in V_i$
  and $\theta(t_i)=t$.  Since $\theta$ sends normal forms to normal
  forms, $\theta$ is injective; and since the image of $\theta$ is
  precisely the set of all normal forms $t^{\pm 1}v_1t^{\pm
    1}v_2\dots$ with $t$-exponent sum divisible by $k$, the image of
  $\theta$ is precisely $K$.

  Next, let $\Delta_*$ be $\Delta_k$ with an additional vertex $*$ and
  $k$ additional edges $(i,*)$ ($i\in\zz/k$).  Taking vertex groups to
  be the $V_i$ and $V_*$ (corresponding to $*$), edge groups to be the
  $E_i$ and copies of the $V_i$ (corresponding to $(i,*)$), and
  attaching homomorphisms to be the $f_i$ and $g_i$, along with the
  identity on $V_i$ (corresponding to the initial vertex of $(i,*)$)
  and $h_i$ (corresponding to the terminal vertex of $(i,*)$), we get
  a graph of groups $(\set{V_*,V_i,E_i},\Delta_*)$.

  Since $(\set{V_i,E_i},\Delta_k)$ is a sub-(graph of groups) of
  $(\set{V_*,V_i,E_i},\Delta_*)$, by the normal form theorem, its
  fundamental group $K$ is a subgroup of the fundamental group 
  $\pi_1(\set{V_*,V_i,E_i},\Delta_*)$.  Furthermore, if $T_\Delta$ is the
  spanning tree of $\Delta_*$ whose edges are the $(i,*)$ and $T_\Gamma$ is
  the (unique) spanning tree of $\Gamma_k$, by
  \eqref{eq:inT}--\eqref{eq:notinT}, the group
  $\pi_1(\set{V_*,V_i,E_i},\Delta_*,T_\Delta)$ is isomorphic to
  $\pi_1(\set{V_*,E_i},\Gamma_k,T_\Gamma)$.  The lemma follows.
\end{proof} 

For the rest of this section, we let $S$, $m$, $\phi: S'\rightarrow
S''$, and $G(S,m,\phi)$ be as defined in Section~\ref{results}.  Let
$\Lambda$ be the (directed) graph with vertex set $S$ and an edge
$(s,\phi(s))$ for each $s\in S'$.  Since every vertex of $\Lambda$ has
in-degree and out-degree each at most 1, $\Lambda$ is a disjoint union
of loops and (non-closed) paths.  Choose an integer $k$ to be a
multiple of the length of each loop in $\Lambda$ and strictly greater
than twice the length of each path in $\Lambda$.

Let $\barS$ be the quotient of the direct product $\zz/k\times S$ by
the equivalence relation generated by $(i,s)\sim(i+1,\phi(s))$ for all
$i\in \zz/k$ and $s\in S'$.  More geometrically, let $\Lambdahat$ be
the (directed) graph with vertex set $\zz/k\times S$ and an edge
$((i,s),(i+1,\phi(s)))$ for each $s\in S'$ and each $i\in \zz/k$.
Then the map $(i,s)\mapsto s$ on vertex sets extends to a $k$-fold
covering map $\pi:\Lambdahat\rightarrow\Lambda$, and each element of
$\barS$ is the vertex set of a connected component of $\Lambdahat$.

For $i\in \zz/k$, define $\psi_i:S\rightarrow \barS$ by letting
$\psi_i(s)$ be the equivalence class of $(i,s)$.  We would like to
define a labelling function $\barm$ on $\barS$ by
\begin{equation}\label{eq:label-function}
  \barm(x,y) =
  \begin{cases}
    m(s,t) &\text{if $x=\psi_i(s)$ and $y=\psi_i(t)$ for some
      $i\in\zz/k$,} \\
    \infty &\text{otherwise.}
  \end{cases}
\end{equation}

\begin{lemma}\label{labels-well-defined}
  The labelling function $\barm$ is well-defined.
\end{lemma}

The proof of Lemma~\ref{labels-well-defined} uses a variation on an
idea from \cite[Lem.\ 5.7]{wise}.

\begin{proof}
  Let $p:(i,s)\mapsto i$ be the natural projection from $\Lambdahat$
  to the graph $\Delta_k$ defined above, and let
  $\pi:\Lambdahat\rightarrow \Lambda$ be the covering map defined
  above.  If $\lhat$ is a component of $\Lambdahat$ such that
  $\ell=\pi(\lhat)$ is a closed loop, then since the length of $\ell$
  divides $k$, $p$ induces a graph isomorphism from $\lhat$ to
  $\Delta_k$.  Similarly, if $\lhat$ is a component of $\Lambdahat$
  such that $\ell=\pi(\lhat)$ is a non-closed path in $\Lambda$, then
  $p$ induces a graph isomorphism from $\lhat$ to $p(\lhat)$, a path
  of length strictly less than $k/2$.  It follows that for any two
  components $\lhat$, $\lhat'$ of $\Lambdahat$, the intersection
  $p(\lhat)\cap p(\lhat')$ is a connected subgraph of $\Delta_k$.

  So now, suppose we have $i,j\in\zz/k$ and $s,s',t,t'\in S$ such that
  $x=\psi_i(s)=\psi_j(s')$ and $y=\psi_i(t)=\psi_j(t')$.
  Equivalently, suppose the vertices $(i,s)$ and $(j,s')$ lie in the
  same component $x$ of $\Lambdahat$ and the vertices $(i,t)$ and
  $(j,t')$ lie in the same component $y$ of $\Lambdahat$.  Since $p$
  induces graph isomorphisms from $x$ and $y$ to their images $p(x)$
  and $p(y)$, respectively, and $p(x)\cap p(y)$ is connected, for some
  $n\ge 0$, either both $(j,s')=(i+n,\phi^n(s))$ and
  $(j,t')=(i+n,\phi^n(t))$, or both $(i,s)=(j+n,\phi^n(s'))$ and
  $(i,t)=(j+n,\phi^n(t'))$.  In either case, by
  \eqref{eq:label-preserving}, we see that $m(s,t)=m(s',t')$.
\end{proof}

\begin{proof}[of Theorem \ref{main}]
  Taking the Artin case first, let $V=A(S,m)$, $E=A(S',m)$, and
  $V_*=A(\barS,\barm)$ (which is well-defined by
  Lemma~\ref{labels-well-defined}).  Also, let $f: E\rightarrow V$ be
  induced by the inclusion of $S'$ in $S$ and let $g: E\rightarrow V$
  be induced by $\phi:S'\rightarrow S''$.  Finally, since $k$ is a
  multiple of the length of each loop in $\Lambda$ and strictly longer
  than each path in $\Lambda$, every lift of a loop or path in
  $\Lambda$ intersects each image $\psi_i(S)$ at most once, which
  means that each of the maps $\psi_i: S\rightarrow\barS$ is injective
  and each $\psi_i$ induces an injection $h_i: V\rightarrow V_*$.
  Applying Lemma~\ref{graphlem} with this choice of $V$, $E$, $V_*$,
  $f$, $g$, and $h_i$, we see that $G(S,m,\phi)$ has a subgroup of
  index $k$ that embeds in $G_*=\pi_1(\set{V_*,E_i},\Gamma_k)$.

  It remains to show that $G_*$ is an Artin group.  By
  \eqref{eq:notinT}, we see that $G_*$ has a presentation with
  generators $\barS\cup\set{t_i}$ ($i\in\zz/k$), subject to the Artin
  relations for $(\barS,\barm)$ together with the relations
  \begin{equation}\label{eq:HNN-is-artin}
    t_i^{-1}\psi_i(s)t_i=\psi_{i+1}(\phi(s))
  \end{equation}
  for all $i\in\zz/k$ and $s\in S'$.  However, since
  $(i,s)=(i+1,\phi(s))$ in $\barS$, \eqref{eq:HNN-is-artin} becomes
  the Artin relation $\psi_i(s)t_i=t_i\psi_i(s)$.  Therefore, if
  \begin{equation}
    m^+(x,y) =
    \begin{cases}
      \barm(x,y) &\text{if $x,y\in\barS$,} \\
      2 &\text{if $\set{x,y}=\set{t_i,\psi_i(s)}$ for some
        $i\in\zz/k$,} \\
      \infty &\text{otherwise,}
    \end{cases}
  \end{equation}
  we see that $G_*=A(\barS\cup\set{t_i},m^+)$.

  Now consider the Coxeter case.  Replacing $V=A(S,m)$, $E=A(S',m)$, and
  $V_*=A(\barS,\barm)$ with $V=W(S,m)$, $E=W(S',m)$, and
  $V_*=W(\barS,\barm)$, the previous proof carries through, except
  that $G_*$ is not a Coxeter group.  However, consider
  $W^+=W(\barS\cup\set{u_i,u'_i},m^+)$, where
  \begin{equation}
    m^+(x,y) =
    \begin{cases}
      \barm(x,y) &\text{if $x,y\in\barS$,} \\
      2 &\text{if $\set{x,y}=\set{u_i,\psi_i(s)}$ for some
        $i\in\zz/k$,} \\
      2 &\text{if $\set{x,y}=\set{u'_i,\psi_i(s)}$ for some
        $i\in\zz/k$,} \\
      \infty &\text{otherwise.}
    \end{cases}
  \end{equation}
  It follows from Tits' solution to the word problem for Coxeter
  groups that the homomorphism $\theta: G_*\rightarrow W^+$ defined by
  $\theta(s)=s$ for $s\in\barS$ and $\theta(t_i)=u_i u'_i$ is
  injective, as reduced (graph of groups) words in $G_*$ are sent to
  reduced (Coxeter) words in $W^+$.  The theorem follows.
\end{proof}

\section{An open question}

Can Theorem~\ref{main} be extended to ``Artin graphs of Artin groups''
(defined similarly to Artin HNN-extensions)?  Note that since an
``Artin amalgamated free product'' is itself an Artin group, it is
enough to consider multiple-Artin-HNN-extensions.  In that case, the
key point is to generalize Lemma~\ref{labels-well-defined}, but it is
not clear to us how (or if) that can be done.

One may of course also consider the analogous Coxeter question(s).

\noindent 
T. Hsu

\noindent
Department of Mathematics
San Jos\'{e} State University
San Jos\'{e}
CA 95192-0103
USA

\noindent 
{\tt hsu@math.sjsu.edu}

\noindent 
I. J. Leary

\noindent 
Department of Mathematics
The Ohio State University
231 W. 18th Avenue
Columbus
OH 43210
USA

\noindent
{\tt leary@math.ohio-state.edu}
\end{document}